\documentclass[leqno,oneside,11pt]{amsart}
\usepackage{amssymb}
\theoremstyle{plain}
	\newtheorem{Thm}{Theorem}[section]
	\newtheorem*{Main}{Main Theorem}
	\newtheorem{Cor}{Corollary}[section]
	\newtheorem{Lem}{Lemma}[section]
	\newtheorem{Prop}{Proposition}[section]
	
	\newtheorem{Claim}{Claim}[section]
\theoremstyle{definition}
	\newtheorem{Def}{Definition}
\theoremstyle{remark}
	 
	\newtheorem{Rem}{Remark}

\numberwithin{equation}{section}
 
\newif\ifShowLabels
\ShowLabelstrue
\newdimen\theight
\def\TeXref#1{%
	\leavevmode\vadjust{\setbox0=\hbox{{\tt
		\quad\quad  {\small \rm #1}}}%
	\theight=\ht0
	\advance\theight by \lineskip
	\kern -\theight \vbox to
	\theight{\rightline{\rlap{\box0}}%
	\vss}%
	}}%

\renewcommand{\sec}[2]{\section{#2}\label{S:#1}%
	\ifShowLabels \TeXref{{S:#1}} \fi}

\newcommand{\refs}[1]{Section ~\ref{S:#1}}

\newcommand{\reft}[1]{Theorem ~\ref{T:#1}}

\newcommand{\refl}[1]{Lemma ~\ref{L:#1}}
\newcommand{\refp}[1]{Proposition ~\ref{P:#1}}
\newcommand{\refc}[1]{Corollary ~\ref{C:#1}}

\newcommand{\refe}[1]{\eqref{E:#1}}

\newenvironment{thm}[1]%
	{ \begin{Thm} \label{T:#1}  \ifShowLabels \TeXref{T:#1} \fi }%
	{ \end{Thm} }
\renewcommand{\th}[1]{\begin{thm}{#1}}
\renewcommand{\eth}{\end{thm} }

	{ \begin{Main} \label{T:main}  \ifShowLabels \TeXref{T:main} \fi }%
	{ \end{Main} }

\newenvironment{lemma}[1]%
	{ \begin{Lem} \label{L:#1}  \ifShowLabels \TeXref{L:#1} \fi }%
	{ \end{Lem} }
\newcommand{\lem}[1]{\begin{lemma}{#1}}
\newcommand{\elem}{\end{lemma}}

\newenvironment{propos}[1]%
	{ \begin{Prop} \label{P:#1}  \ifShowLabels \TeXref{P:#1} \fi }%
	{ \end{Prop} }
\newcommand{\prop}[1]{\begin{propos}{#1}}
\newcommand{\eprop}{\end{propos}}

\newenvironment{corol}[1]%
	{ \begin{Cor} \label{C:#1}  \ifShowLabels \TeXref{C:#1} \fi }%
	{ \end{Cor} }
\newcommand{\cor}[1]{\begin{corol}{#1}}
\newcommand{\ecor}{\end{corol}}

\newenvironment{defeni}[1]%
	{ \begin{Def} \label{D:#1}  \ifShowLabels \TeXref{D:#1} \fi }%
	{ \end{Def} }
\newcommand{\defe}[1]{\begin{defeni}{#1}}
\newcommand{\edefe}{\end{defeni}}

\newenvironment{claimenv}[1]%
	{ \begin{Claim} \label{C:#1}  \ifShowLabels \TeXref{C:#1} \fi }%
	{ \end{Claim} }
\newcommand{\claim}[1]{\begin{claimenv}{#1}}
\newcommand{\eclaim}{\end{claimenv}}

\newenvironment{remark}[1]%
	{ \begin{Rem} \label{R:#1}  \ifShowLabels \TeXref{R:#1} \fi }%
	{ \end{Rem} }
\newcommand{\rem}[1]{\begin{remark}{#1}}
\newcommand{\erem}{\end{remark}}

	{ \begin{equation} \label{E:#1}  \ifShowLabels \TeXref{E:#1} \fi }%
	{ \end{equation} }
\newcommand{\eq}[1]%
	{ \ifShowLabels \TeXref{E:#1} \fi
	  \begin{equation} \label{E:#1} }
\newcommand{\eeq}{\end{equation}}

\newcommand{\prf}{ \begin{proof} }
\newcommand{\eprf}{ \end{proof} }

\newcommand\alp{\alpha}

\newcommand\bet{\beta}

\newcommand\Gam{\Gamma}
\newcommand\del{\delta}

\newcommand\tet{\theta}

\newcommand\lam{\lambda}

\newcommand\calI{{\mathcal{I}}}

\newcommand\calL{{\mathcal{L}}}

\newcommand\calO{{\mathcal{O}}}

\newcommand\bfOme{{\mathbf \Omega}}

\newcommand\CC{\mathbb{C}}

\newcommand\fru{{\mathfrak{u}}}

\newcommand\rank{\operatorname{rank}}
\newcommand\irr{\operatorname{Irr}}
\renewcommand{\hom}{\operatorname{Hom}}
\newcommand{\cf}{\operatorname{cf}}
\newcommand{\rad}{\operatorname{Rad}}

\newcommand{\lset}{\left\{}
\newcommand{\rset}{\right\}}
\newcommand{\lpar}{\left(}
\newcommand{\rpar}{\right)}
\newcommand{\lbra}{\left[}
\newcommand{\rbra}{\right]}
\newcommand{\lan}{\left\langle}
\newcommand{\ran}{\right\rangle}
\renewcommand{\lvert}{\left|}
\renewcommand{\rvert}{\right|}
\newcommand{\rar}{\rightarrow}
\newcommand{\fq}{\mathbb{F}_{q}}
\newcommand{\inv}{^{-1}}
\newcommand{\x}{\times}
\newcommand{\+}{\oplus}
\newcommand{\bas}[2]{\lpar {#1}_{1}, \ldots, {#1}_{#2} \rpar}
\newcommand{\seq}[2]{{#1}_{1}, \ldots, {#1}_{#2}}
\newcommand{\all}{\text{ for all }}
\newcommand{\pr}{^{\prime}}
\newcommand{\sset}{\subseteq}
\newcommand\dsty[1]{{\displaystyle #1}}
\newcommand{\bca}{\begin{cases}}
\newcommand{\eca}{\end{cases}}
\newcommand{\al}[1]{\item[{\rm (#1)}]}
\newcommand{\ie}{i.e., }
\newcommand{\mat}{\begin{bmatrix}}
\newcommand{\emat}{\end{bmatrix}}
\newcommand{\smat}{\left[ \begin{smallmatrix}}
\newcommand{\esmat}{\end{smallmatrix} \right]}
\newcommand{\bcup}{\bigcup}
\newcommand{\col}{\colon}
\newcommand{\orb}{\calO}
\newcommand{\Orb}{\bfOme}
\newcommand{\Sset}{\supseteq}
\newcommand{\ord}[1]{\lvert #1 \rvert}
\newcommand{\ovl}{\overline}
\newcommand{\imp}[2]{\text{(#1)} \Rightarrow \text{(#2)}}
\newcommand{\impd}[2]{\text{(#1)} \Leftrightarrow \text{(#2)}}

\newcommand{\apl}[3]{#1 \col #2 \rar #3}

\ShowLabelsfalse

\begin{document}

\title[Irreducible characters of finite algebra groups]{Irreducible 
characters of finite algebra groups}

\author[C.~A.~M.~Andr\'e]{Carlos A. M. Andr\'e}

\thanks{This research was carried out as part of the PRAXIS XXI 
Project 2/2.1/MAT/73/94.}

\address{Departamento de Matem\'atica, Faculdade de Ci\^encias da 
Universidade de Lisboa, Rua Ernesto de Vasconcelos, Bloco C1, Piso 3, 
1700 LISBOA, PORTUGAL}

\email{candre@fc.ul.pt}

\maketitle

\sec{intro}{Introduction}

Let $p$ be a prime number, let $q = p^{e}$ ($e \geq 1$) be a power of 
$p$ and let $\fq$ denote the finite field with $q$ elements. Let $A$ be 
a finite dimensional $\fq$-algebra. (Throughout the paper, all algebras 
are supposed to have an identity element). Let $J = J(A)$ be the Jacobson 
radical of $A$ and let $$G = 1+J = \lset 1+a \col a \in A \rset.$$ 
Then $G$ is a $p$-subgroup of the group of units of $A$. Following 
\cite{isaacs1}, we refer to a group arising in this way as an {\bf
$\fq$-algebra group}. As an example, let $J = \fru_{n}(q)$ be the 
$\fq$-space consisting of all nilpotent uppertriangular $n \x n$ 
matrices over $\fq$. Then $J$ is the Jacobson radical of the 
$\fq$-algebra $A = \fq \cdot 1 + J$ and the $p$-group $G = 1+J$ is the 
group $U_{n}(q)$ consisting of all unipotent uppertriangular $n 
\x n$ matrices over $\fq$.

A subgroup $H$ of an $\fq$-algebra group $G$ is said to be an {\bf 
algebra subgroup} of $G$ if $H = 1+U$ for some multiplicatively closed 
$\fq$-subspace $U$ of $J$. It is clear that an algebra subgroup of $G$ 
is itself an $\fq$-algebra group and that it has $q$-power index in $G$. 

The main purpose of this paper is to proof the following result. 
(Throughout this paper, all characters are taken over the complex 
field.)

\th{main}
	Let $G$ be an $\fq$-algebra group and let $\chi$ be an irreducible 
	character of $G$. Then there exist an algebra subgroup $H$ of $G$ and 
	a linear character $\lam$ of $H$ such that $\chi = \lam^{G}$.
\eth

As a consequence, we obtain Theorem~A of \cite{isaacs1} (see also 
\cite[Theorem[26.7]{huppert}) which asserts 
that all irreducible characters of an (arbitrary) $\fq$-algebra group 
have $q$-power degree. (However, this result will used in the proof of 
\reft{main}.) Following the terminology of \cite{isaacs1}, we say that 
a finite group $G$ is a {\it $q$-power-degree group} if every 
irreducible character of $G$ has $q$-power degree. Hence, 
\cite[Theorem~A]{isaacs1} asserts that every $\fq$-algebra group is a 
$q$-power-degree group. In particular, the unitriangular group $U_{n}(q)$ 
is a $q$-power-degree group (which is precisely the statement of 
\cite[Corollary~B]{isaacs1}). On the other hand, our \reft{main} generalizes 
Theorem~C of \cite{isaacs1} and answers the question made by I.~M.~Isaacs 
immediately before that theorem. We note, moreover, that the statement 
of our \reft{main} is precisely the assertion made by E.~A.~Gutkin in 
\cite{gutkin}. The argument used by Gutkin to prove this assertion was 
defective and a counterexample was given by Isaacs to illustrate its flaw 
(see Section~10 of \cite{isaacs1}).

A result similar to our \reft{main} was proved by D.~Kazhdan for the 
group $G = U_{n}(q)$ in the case where $p \geq n$. Kazhdan's result appears 
in the paper \cite{kazhdan} (see also \cite[Theorem~7.7]{srinivasan}) and 
applies to other finite unipotent algebraic groups. However, Kazhdan imposes 
a restriction on the prime $p$ in order to use the exponential map.
In this paper, we replace the exponential map by the bijection $J 
\rar 1+J$ defined by the (natural) correspondence $a \mapsto 
1+a$. Then we follow Kazhdan's idea and we use Kirillov's method of 
coadjoint orbits (see, for example, \cite{kirillov1}) to parametrize 
the irreducible characters of the $\fq$-algebra group $G = 1+J$. 

\sec{class}{Class functions associated with coadjoint orbits}

Let $J = J(A)$, where $A$ is a finite dimensional $\fq$-algebra, and 
let $G = 1+J$. Let $J^{*} = \hom_{\fq}(J, \fq)$ be the dual space of 
$J$ and let $\psi$ be an arbitrary non-trivial linear character of 
the additive group ${\fq}^{+}$ of the field $\fq$. For each $f \in 
J^{*}$, let $\psi_{f} \col J \rar \CC$ be the map defined by 
\eq{psi}
\psi_{f}(a) = \psi(f(a))
\end{equation}
for all $a \in J$. Then $\psi_{f}$ is a linear character of the additive 
group $J^{+}$ of $J$ and, in fact,
\eq{irrJ}
\irr(J^{+}) = \lset \psi_{f} \col f \in J^{*} \rset.
\end{equation}
(For any finite group $X$, we denote by $\irr(X)$ the set of all 
irreducible characters of $X$.)

The group $G$ acts on $J^{*}$ by $(x \cdot f)(a) = f(x\inv a x)$ for 
all $x \in G$, all $f \in J^{*}$ and all $a \in J$. (Usually, we refer 
to this action as the {\bf coadjoint action} of $G$.) Let $\Orb(G)$ 
denote the set of all $G$-orbits on $J^{*}$. We claim that the 
cardinality $\ord{\orb}$ of any $G$-orbit $\orb \in \Orb(G)$ is a 
$q^{2}$-power. To see this, let $f \in J^{*}$ be arbitrary and 
define $B_{f} \col J \x J \rar \fq$ by
\eq{Bf}
B_{f}(a,b) = f([ab])
\end{equation}
for all $a,b \in J$ (here $[ab]=ab-ba$ is the usual Lie product of $a,b 
\in J$). Then $B_{f}$ is a skew-symmetric $\fq$-bilinear form. Let $n = 
\dim J$, let $\bas{e}{n}$ be an $\fq$-basis of $J$ and let $M(f)$ be the 
skew-symmetric matrix which represents $B_{f}$ with respect to this basis. 
Then $M(f)$ has even rank (see, for example, \cite[Theorem~8.6.1]{cohn}). 
Let $$\rad(f) = \lset a \in J \col f([ab]) = 0 \all b \in J \rset$$ be 
the radical of $B_{f}$. Then $\rad(f)$ is an $\fq$-subspace of $J$ and
\eq{rad}
\dim \rad(f) = \dim J - \rank M(f).
\end{equation}
We have the following result.

\prop{rad}
	Let $f \in J^{*}$ be arbitrary. Then $\rad(f)$ is a multiplicatively 
	closed $\fq$-subspace of $J$. Moreover, the centralizer $C_{G}(f)$ 
	of $f$ in $G$ is the algebra subgroup $1+\rad(f)$ of $G$. In 
	particular, if $\orb \in \Orb(G)$ is the $G$-orbit which contains 
	$f$, then $\ord{\orb} = q^{\rank M(f)}$ is a $q^{2}$-power.
\eprop

\begin{proof}
	Since $[ab,c] = [a,bc] + [b,ca]$, we clearly have $f([ab,c]) = 0$ 
	for all $a, b \in \rad(f)$ and all $c \in J$. Thus $\rad(f)$ is 
	multiplicatively closed.
	
	On the other hand, let $x \in G$ be arbitrary. Then $x \in C_{G}(f)$ 
	if and only if $f(x\inv bx) = f(b)$ for all $b \in J$. Hence $x \in 
	C_{G}(f)$ if and only if $f(bx) = f(xb)$ for all $b \in J$. Now, let 
	$a = x-1 \in J$. Then $f(bx) = f(b) + f(ba)$ and $f(xb) = f(b) + 
	f(ab)$, and so $f([ab]) = f(xb) - f(bx)$ for all $b \in J$. It 
	follows that $x \in C_{G}(f)$ if and only if $a \in \rad(f)$.
	
	For the last assertion, we note that $\ord{G} = \ord{C_{G}(f)} \cdot 
	\ord{\orb}$, that $\ord{G} = q^{\dim J}$ and (as we have just proved) 
	that $\ord{C_{G}(f)} = q^{\dim \rad(f)}$. Therefore, by \refe{rad}, 
	we deduce that $\ord{\orb} = q^{\dim J - \dim \rad(f)} = q^{\rank M(f)}$.
\end{proof}

For each $\orb \in \Orb(G)$, we define the function $\phi_{\orb} \col 
G \rar \CC$ by the rule
\eq{class}
	\phi_{\orb}(1+a) = \frac{1}{\sqrt{\ord{\orb}}} \sum_{f \in \orb} 
	\psi_{f}(a)
\end{equation}
for all $a \in J$. It is clear that $\phi_{\orb}$ is a class function of 
$G$ of degree
\eq{degree}
\phi_{\orb}(1) = \sqrt{\ord{\orb}} = q^{\rank M(f)}
\end{equation}
where $M(f)$ is as before. Moreover, we have the following result.

\prop{ortho}
	The set $\lset \phi_{\orb} \col \orb \in \Orb(G) \rset$ is an orthonormal 
	basis for the $\CC$-space $\cf(G)$ consisting of all class functions on 
	$G$. In particular, we have $$\frac{1}{\ord{G}} \sum_{x \in G} 
	\phi_{\orb}(x) \ovl{\phi_{\orb\pr}(x)} = \del_{\orb, \orb\pr}$$ for all 
	$\orb, \orb\pr \in \Orb(G)$. (Here, $\del$ denotes the usual Kronecker 
	symbol.)
\eprop

\begin{proof}
	Let $\lan \cdot , \cdot \ran_{G}$ denote the Frobenius scalar product 
	on $\cf(G)$. Let $\orb, \orb\pr \in \Orb(G)$ be arbitrary. Then
	\begin{eqnarray*}
		\lan \phi_{\orb}, \phi_{\orb\pr} \ran_{G} & = & \frac{1}{\ord{G}} 
		\sum_{x \in G} \phi_{\orb}(x) \ovl{\phi_{\orb\pr}(x)}  \\
		& = & \frac{1}{\ord{J}}\sum_{a \in J} \frac{1}{\sqrt{\ord{\orb}}} 
		\frac{1}{\sqrt{\ord{\orb\pr}}} \sum_{f \in \orb} \sum_{f\pr \in \orb\pr} 
		\psi_{f}(a) \ovl{\psi_{f\pr}(a)}  \\
		& = & \frac{1}{\sqrt{\ord{\orb}}} \frac{1}{\sqrt{\ord{\orb\pr}}} 
		\sum_{f \in \orb} \sum_{f\pr \in \orb\pr} \lpar 
		\frac{1}{\ord{J}}\sum_{a \in J} \psi_{f}(a) \ovl{\psi_{f\pr}(a)} 
		\rpar  \\
	\end{eqnarray*}
	\begin{eqnarray*}
		\phantom{\lan \phi_{\orb}, \phi_{\orb\pr} \ran_{G}}
		& = & \frac{1}{\sqrt{\ord{\orb}}} \frac{1}{\sqrt{\ord{\orb\pr}}} 
		\sum_{f \in \orb} \sum_{f\pr \in \orb\pr} \lan \psi_{f}, \psi_{f\pr} 
		\ran_{J^{+}}  \\
		& = & \frac{1}{\sqrt{\ord{\orb}}} \frac{1}{\sqrt{\ord{\orb\pr}}} 
		\sum_{f \in \orb} 
		\sum_{f\pr \in \orb\pr} \del_{f,f\pr}
	\end{eqnarray*}
	(using \refe{irrJ}) and so $$\lan \phi_{\orb}, \phi_{\orb\pr} \ran_{G} = 
	\del_{\orb, \orb\pr}.$$
	
	To conclude the proof, we claim that $\ord{\Orb(G)}$ equals the class 
	number $k_{G}$ of $G$; we recall that $k_{G} = \dim_{\CC} \cf(G)$ 
	(see, for example, \cite[Corollary~2.7 and Theorem~2.8]{isaacs2}). 
	Firstly, we observe that $k_{G}$ is the number of $G$-orbits on $J$ for 
	the {\it adjoint action}: $x \cdot a = xax\inv$ for all $x \in G$ and 
	all $a \in J$. Let $\tet$ be the permutation character of $G$ on $J$ 
	(see \cite{isaacs2} for the definition). Then, by 
	\cite[Corollary~5.15]{isaacs2}, $$k_{G} = \lan \tet, 1_{G} \ran_{G}.$$ 
	Moreover, by definition, we have $$\tet(x) = \ord{\lset a \in J \col x 
	\cdot a = a \rset}$$ for all $x \in G$.
	
	On the other hand, consider the action of $G$ on $\irr(J^{+})$ given by 
	$$x \cdot \psi_{f} = \psi_{x \cdot f}$$ for all $x \in G$ and all $f \in 
	J^{*}$. We clearly have $$(x \cdot \psi_{f})(x \cdot a) = \psi_{f}(a)$$ 
	for all $x \in G$, all $f \in J^{*}$ and all $a \in J$. It follows from 
	Brauer's Theorem (see \cite[Theorem~6.32]{isaacs2}) that $$\tet(x) = 
	\ord{\lset f \in J^{*} \col x \cdot \psi_{f} = \psi_{f} \rset}$$ for all 
	$x \in G$. Therefore, $\tet$ is also the permutation character of $G$ 
	on $\irr(J^{+})$ and so $$\lan \tet, 1_{G} \ran_{G} = \ord{\Orb(G)}.$$ 
	
	The claim follows and the proof is complete.
\end{proof}

We will prove (see \reft{irred} in \refs{irred}) that $$\irr(G) = \lset 
\phi_{\orb} \col \orb \in \Orb(G) \rset.$$ (This is the key for the proof 
of \reft{main}.) Therefore, the next result will of course be a 
consequence of that theorem. However, we give below a very easy proof 
(independent of \reft{irred}) of the second orthogonality relations for 
the functions $\phi_{\orb}$ for $\orb \in \Orb(G)$.

\prop{ortho2}
	Let $x, y \in G$ be arbitrary. Then $$\sum_{\orb \in \Orb(G)} 
	\phi_{\orb}(x) \ovl{\phi_{\orb}(y)} = \bca \ord{C_{G}(x)}, & 
	\text{if $x$ and $y$ are $G$-conjugate,} \\ 0, & \text{otherwise.} 
	\eca$$
\eprop

\begin{proof}
	Let $a = x - 1$ and $b = y - 1$. Then
	\begin{eqnarray*}
		\sum_{\orb \in \Orb(G)} \phi_{\orb}(x) \ovl{\phi_{\orb}(y)} & = & 
		\sum_{\orb \in \Orb(G)} \frac{1}{\ord{\orb}} \sum_{f \in \orb} \sum_{g 
		\in \orb} \psi_{f}(a) \ovl{\psi_{g}(b)}  \\
		& = & \sum_{\orb \in \Orb(G)} \frac{1}{\ord{\orb}} \sum_{f \in \orb} 
		\frac{1}{\ord{C_{G}(f)}} \sum_{z \in G} \psi_{f}(a) 
		\ovl{\psi_{f}(z\inv b z)}  \\
		& = & \sum_{\orb \in \Orb(G)} \sum_{f \in \orb} \frac{1}{\ord{G}} 
		\sum_{z \in G} \psi_{f}(a - z\inv bz)  \\
		& = & \frac{1}{\ord{G}} \sum_{z \in G} \sum_{f \in J^{*}} 
		\psi_{f}(a - z\inv bz)  \\
		& = & \frac{1}{\ord{G}} \sum_{z \in G} \rho_{J^{+}}(a - z\inv bz)
	\end{eqnarray*}
	where $\rho_{J^{+}}$ denotes the regular character of $J^{+}$. 
	It follows that $$\sum_{\orb \in \Orb(G)} \phi_{\orb}(x) 
	\ovl{\phi_{\orb}(y)} = \sum_{z \in G} \del_{a,z\inv bz} = \ord{\lset z 
	\in G \col a = z\inv bz \rset}$$ and this clearly completes the proof.
\end{proof}

As a consequence we obtain the following additive decomposition of the 
regular character $\rho_{G}$ of $G$ (which is also a consequence of 
\reft{irred}).

\cor{regular}
	$\dsty{\rho_{G} = \sum_{\orb \in \Orb(G)} \phi_{\orb}(1) \phi_{\orb}}$.
\ecor

\begin{proof}
	Let $x \in G$ be arbitrary. By the previous proposition, 
	$$\sum_{\orb \in \Orb(G)} \phi_{\orb}(1) \phi_{\orb}(x) = 
	\del_{x,1} \ord{G} = \rho_{G}(x).$$ The result follows (by 
	definition of $\rho_{G}$).
\end{proof}

\sec{max}{Maximal algebra subgroups}

In this section, we consider restriction and induction of the class 
functions defined in the previous section. We follow Kirillov's 
theory on nilpotent Lie groups (see, for example, 
\cite{corwin-greenleaf}). As before, let $A$ be a finite 
dimensional $\fq$-algebra, let $J = J(A)$ and let $G = 1+J$.

Let $U$ be a maximal multiplicatively closed $\fq$-subspace of $J$. 
Then $J^{2} \sset U$; otherwise, we must have $U + J^{2} = J$ and this 
implies that $U = J$ (see \cite[Lemma~3.1]{isaacs1}). It follows that 
$U$ is an ideal of $A$ and so $H = 1+U$ is a normal subgroup of $G$ (in 
the terminology of \cite{isaacs1}, we say that $H$ is an {\bf ideal 
subgroup} of $G$). Moreover, we have $\dim U = \dim J - 1$ and so 
$\ord{G:H} = q$.

Let $\pi \col J^{*} \rar U^{*}$ be the natural projection (by definition, 
for any $f \in J^{*}$, $\pi(f) \in U^{*}$ is the restriction of $f$ to $U$). 
Then the kernel of $\pi$ is the $\fq$-subspace $$U^{\perp} = \lset f \in 
J^{*} \col f(a) = 0 \all a \in U \rset$$ of $U$. On the other hand, for 
any $f \in U^{*}$, the fibre $\pi\inv(\pi(f))$ of $\pi(f) \in U^{*}$ is 
the subset $$\calL(f) = \lset g \in U^{*} \col g(a) = f(a) \all a \in U 
\rset$$ of $J^{*}$. It is clear that $$\calL(f) = f + U^{\perp} = \lset 
f+g \col g \in U^{\perp} \rset$$ for all $f \in J^{*}$.

Let $f \in J^{*}$ be arbitrary and let $f_{0}$ denote the projection 
$\pi(f) \in U^{*}$. Let $\orb \in \Orb(G)$ be the $G$-orbit which 
contains $f$ and let $\orb_{0} \in \Orb(H)$ be the $H$-orbit which 
contains $f_{0}$. Since $\pi(x \cdot f) = x \cdot \pi(f)$ for all $x 
\in G$ (because $H$ is normal in $G$, hence $U$ is invariant under the 
adjoint $G$-action), the projection $\pi(\orb) \sset U^{*}$ of $\orb$ 
is $G$-invariant. Thus $\orb_{0} \sset \pi(\orb)$; in fact, $\pi(\orb)$ 
is a disjoint union of $H$-orbits. It follows $$\ord{\orb_{0}} \leq 
\ord{\pi(\orb)} \leq \ord{\orb}.$$ Let $\pi_{\orb} \col \orb \rar 
\pi(\orb)$ denote the restriction of $\pi$ to $\orb$. Since $\pi_{\orb}$ 
is surjective, $\orb$ is the disjoint union $$\orb = \bigcup_{g_{0} \in 
\pi(\orb)} \pi\inv(g_{0}).$$ Since $\pi\inv(\pi(g)) = \calL(g) \cap 
\orb$ for all $g \in \orb$, we conclude that $$\ord{\orb} = 
\sum_{g \in \Gam_{\orb}} \ord{\calL(g) \cap \orb},$$ where $\Gam_{\orb} 
\sset \orb$ is a set of representatives of the fibres $\pi\inv(g_{0})$ 
for $g_{0} \in \pi(\orb)$. Now, we clearly have $\calL(x \cdot g) = 
x \cdot \calL(g)$ for all $x \in G$ and all $g \in J^{*}$. It follows 
that
\eq{inters}
	x \cdot \lpar \calL(g) \cap \orb \rpar = \calL(x \cdot g) \cap 
	\orb
\end{equation}
for all $x \in G$ and all $g \in J^{*}$. In particular, we have 
$\ord{\calL(g) \cap \orb} = \ord{\calL(f) \cap \orb}$ for all $g \in \orb$. 
Since $\ord{\Gam_{\orb}} = \ord{\orb}$, we conclude that
\eq{order}
	\ord{\orb} = \ord{\pi(\orb)} \cdot \ord{\calL(f) \cap \orb}.
\end{equation}
In particular, we deduce that $$\ord{\orb} \leq q \ord{\pi(\orb)}$$ 
(because  $\calL(f) = f + U^{\perp}$, hence $\ord{\calL(f)} = q$). We 
claim that $\ord{\pi(\orb)}$ is a power of $q$. In fact, we have the 
following.

\lem{rad}
	Let $H$ be a maximal algebra subgroup of $G$, let $U \sset J$ be such 
	that $H = 1+U$ and let $\apl{\pi}{J^{*}}{U^{*}}$ be the natural 
	projection. Let $f \in J^{*}$ be arbitrary and let $\orb \in \Orb(G)$ 
	be the $G$-orbit which contains $f$. Then:
	\begin{enumerate}
		\al{a}  The centralizer $C_{G}(f_{0})$ of $f_{0} = \pi(f)$ in $G$ is 
		an algebra subgroup of $G$; in fact, $C_{G}(f_{0}) = 1 + R$ where $R = 
		\lset a \in J \col f([ab]) = 0 \all b \in U \rset$.
		\al{b}  $\ord{\pi(\orb)}$ is a power of $q$; in fact, either 
		$\ord{\orb} = \ord{\pi(\orb)}$, or $\ord{\orb} = q \ord{\pi(\orb)}$.
	\end{enumerate}
\elem

\begin{proof}
	The proof of (a) is analoguous to the proof of \refp{rad}. On the 
	other hand, since $G$ acts transitively on $\pi(\orb)$, we have 
	$$\ord{\pi(\orb)} = \ord{G} \cdot \ord{C_{G}(f_{0})}\inv = q^{\dim J - 
	\dim R}.$$ The last assertion is clear because $\ord{\pi(\orb)} 
	\leq \ord{\orb} \leq q \ord{\pi(\orb)}$.
\end{proof}

Following \cite{kirillov2}, we say a $G$-orbit $\orb \in \Orb(G)$ is of 
{\bf type I} (with respect to $H$) if $\ord{\orb} = \ord{\pi(\orb)}$; 
otherwise, if $\ord{\orb} = q \ord{\pi(\orb)}$, we say that $\orb$ is of 
{\bf type II} (with respect to $H$). The following result asserts that 
our definition is, in fact, equivalent to Kirillov's definition.

\prop{types}
	Let $H$ be a maximal algebra subgroup of $G$, let $U \sset J$ be such 
	that $H = 1+U$ and let $\apl{\pi}{J^{*}}{U^{*}}$ be the natural projection. 
	Let $\orb \in \Orb(G)$ be arbitrary. Then:
	\begin{enumerate}
		\al{a}  $\orb$ is of type I (with respect to $H$) if and only 
		if $\calL(f) \cap \orb = \lset f \rset$ for all $f \in \orb$;
		\al{b}  $\orb$ is of type II (with respect to $H$) if and only 
		if $\calL(f) \sset \orb$ for all $f \in \orb$.
	\end{enumerate}
\eprop

\begin{proof}
	Suppose that $\ord{\orb} = \ord{\pi(\orb)}$ (\ie $\orb$ is of type I). 
	Then \refe{order} implies that $\ord{\calL(f) \cap \orb} = 1$ 
	and so $\calL(f) \cap \orb = \lset f \rset$. On the other hand, 
	if $\ord{\orb} = q \ord{\pi(\orb)}$ (\ie if $\orb$ is of type II), we 
	must have $\ord{\calL(f) \cap \orb} = q$ and so $\calL(f) \sset \orb$ 
	(because $\ord{\calL(f)} = q$). The result follows by 
	\refe{inters}.
\end{proof}

Now, let $n = \dim J$ and let $\bas{e}{n}$ be an $\fq$-basis of $J$ such 
that $e_{i} \in U$ for all $1 \leq i \leq n-1$. Moreover, let $M = 
M(f)$ be the $n \x n$ skew-symmetric matrix which represents the bilinear 
form $B_{f}$ with respect to the basis $\bas{e}{n}$. By \refp{rad}, 
$\ord{\orb} = q^{\rank M}$. Moreover, the matrix $M$ has the form $$M =
\mat M_{0} & - \nu^{T} \\ \nu & 0 \emat$$ where $M_{0} = M(f_{0})$ is the 
$(n-1) \x (n-1)$ skew-symmetric matrix which represents the bilinear form 
$B_{f_{0}} \col U \times U \rar \fq$ with respect to the $\fq$-basis 
$\bas{e}{n-1}$ of $U$, and $\nu$ is the row vector $\nu = \lbra 
f([e_{n}\, e_{1}])\ \cdots\ f([e_{n}\, e_{n-1}]) \rbra$. Since 
$\orb_{0}$ is the $H$-orbit of the element $f_{0} \in U^{*}$, we 
have $\ord{\orb_{0}} = q^{\rank M_{0}}$ (by \refp{rad}). Since $M$ and 
$M_0$ are skew-symmetric matrices, they have even ranks and so, either 
$\rank M = \rank M_{0}$, or $\rank M = \rank M_{0} + 2$. This concludes 
the proof of the following.

\lem{order}
	Let $H$ be a maximal algebra subgroup of $G$, let $U \sset J$ be such 
	that $H = 1+U$ and let $\apl{\pi}{J^{*}}{U^{*}}$ be the natural 
	projection. Let $\orb \in \Orb(G)$ be arbitrary and let $\orb_{0} 
	\in \Orb(H)$ be such that $\orb_{0} \sset \pi(\orb)$. Then, either 
	$\ord{\orb} = \ord{\orb_{0}}$, or $\ord{\orb} = q^{2} \ord{\orb_{0}}$.
\elem

We note that, since $\ord{\orb_{0}} \leq \ord{\pi(\orb)} \leq \ord{\orb}$, 
the equality $\ord{\orb} = \ord{\orb_{0}}$ implies that $\orb$ is of type 
I (with respect to $H$); hence, $\ord{\orb} = q^{2} \ord{\orb_{0}}$ 
whenever $\orb$ is of type II (with respect to $H$). Our next result shows 
that the dicothomy of the preceeding lemma characterizes the $G$-orbit 
$\orb$ with respect to the subgroup $H$.

\prop{char}
	Let $H$ be a maximal algebra subgroup of $G$, let $U \sset J$ be such 
	that $H = 1+U$ and let $\apl{\pi}{J^{*}}{U^{*}}$ be the natural 
	projection. Let $\orb \in \Orb(G)$ be arbitrary and let $\orb_{0} 
	\in \Orb(H)$ be such that $\orb_{0} \sset \pi(\orb)$. Moreover, let $f 
	\in \orb$ be such that the projection $f_{0} = \pi(f)$ lies in 
	$\orb_{0}$. Then:
	\begin{enumerate}
		\al{a}  The following are equivalent:
		\begin{enumerate}
			\al{i}  $\ord{\orb}$ is of type I (with respect to $H$);
			\al{ii}  $\ord{\orb} = \ord{\orb_{0}}$;
			\al{iii}  $\dim \rad(f) = \dim \rad(f_{0}) + 1$;
			\al{iv}  $\ord{C_{G}(f)} = q \ord{C_{H}(f_{0})}$.
		\end{enumerate}
		\al{b}  The following are equivalent:
		\begin{enumerate}
			\al{i}  $\ord{\orb}$ is of type II (with respect to $H$);
			\al{ii}  $\ord{\orb} = q^{2} \ord{\orb_{0}}$;
			\al{iii}  $\dim \rad(f) = \dim \rad(f_{0}) - 1$;
			\al{iv}  $\ord{C_{G}(f)} = q\inv \ord{C_{H}(f_{0})}$.
		\end{enumerate}
	\end{enumerate}
\eprop

\begin{proof}
	The equivalence $\impd{iii}{iv}$ (in both (a) and (b)) follows from 
	\refp{rad}. On the other hand, the equivalence $\impd{ii}{iii}$ (in 
	both (a) and (b)) follows from \refe{rad} (using also \refp{rad}). 
	We have already proved that $\imp{ii}{i}$ in (a) (which is equivalent 
	to $\imp{i}{ii}$ in (b)). Conversely, suppose that $\ord{\orb}$ is of 
	the type I (with respect to $H$). Then $\ord{\pi(\orb)} = \ord{\orb}$. 
	Since $G$ acts transitively on $\pi(\orb)$, we deduce that $$\ord{G} = 
	\ord{C_{G}(f_{0})} \cdot \ord{\pi(\orb)} = \ord{C_{G}(f_{0})} \cdot 
	\ord{\orb}.$$ It follows that $\ord{C_{G}(f_{0})} = \ord{C_{G}(f)}$. 
	Since $C_{G}(f) \sset C_{G}(f_{0})$, we conclude that $C_{G}(f_{0}) = 
	C_{G}(f)$. On the other hand, $C_{H}(f_{0}) \sset C_{G}(f_{0})$ and so 
	$C_{H}(f_{0}) \sset C_{G}(f)$. By the equivalence $\impd{ii}{iv}$ (in 
	both (a) and (b)), we conclude that $\ord{\orb} = \ord{\orb_{0}}$. 
	This completes the proof of $\imp{i}{ii}$ in (a). Hence, the 
	implication $\imp{ii}{i}$ in (b) is also true. The proof is complete.
\end{proof}

We note that, since $\orb_{0} \sset \pi(\orb)$, the equality $\ord{\orb} = 
\ord{\orb_{0}}$ implies that $\pi(\orb) = \orb_{0}$ (hence, the 
projection $\pi(\orb)$ of the $G$-orbit $\orb$ is an $H$-orbit). This 
concludes the proof of part (a) of the following result.

\prop{dec}
	Let $H$ be a maximal algebra subgroup of $G$, let $U \sset J$ be such 
	that $H = 1+U$ and let $\apl{\pi}{J^{*}}{U^{*}}$ be the natural 
	projection. Let $\orb \in \Orb(G)$ be arbitrary and let $\orb_{0} 
	\in \Orb(H)$ be such that $\orb_{0} \sset \pi(\orb)$. Moreover, let 
	$f \in \orb$ be such that the projection $f_{0} = \pi(f)$ lies in 
	$\orb_{0}$. Then the following statements hold:
	\begin{enumerate}
		\al{a}  If $\orb$ is of type I (with respect to $H$), then 
		$\pi(\orb) = \orb_{0}$ is a single $H$-orbit on $U^{*}$.
		\al{b}  Suppose that $\orb$ is of type II (with respect to $H$). Let 
		$e \in J$ be such that $J = U \+ \fq e$ and, for each $\alp \in 
		\fq$, let $x_{\alp}$ denote the element $1 + \alp e \in G$. Then 
		$\pi(\orb)$ is the disjoint union $$\pi({\orb}) = 	\bcup_{\alp \in 
		\fq} \orb_{\alp}$$ where, for each $\alp \in \fq$, $\orb_{\alp} \sset 
		U^{*}$ is the $H$-orbit which contains the element $x_{\alp} \cdot 
		f_{0} \in U^{*}$. We have $$\orb_{\alp} = x_{\alp} \cdot \orb_{0} = 
		\lset x_{\alp} \cdot g \col g \in \orb_{0} \rset$$ for all $\alp \in 
		\fq$. Moreover, the set $\lset x_{\alp} \col \alp \in \fq \rset$ 
		can be replaced by any set of representatives of the of cosets of 
		$H$ in $G$.
	\end{enumerate}
\eprop

\begin{proof}
	It remains to prove part (b).
	
	Let $\alp \in \fq$ be arbitrary. Since $\pi$ is $G$-invariant, we have 
	$x_{\alp} \cdot \pi(f) = \pi(x_{\alp} \cdot f)$ and so $x_{\alp} \cdot 
	f_{0} \in \pi(\orb)$ (we recall that $f_{0} = \pi(f)$). It follows 
	that $\orb_{\alp} \sset \pi(\orb)$.
	
	Next, we prove that $\orb_{\alp} = x_{\alp} \cdot \orb_{0}$. To see this, 
	let $x \in H$ be arbitrary. Then $x \cdot (x_{\alp} \cdot f_{0}) = (x 
	x_{\alp}) \cdot f_{0} =  (x_{\alp} x_{\alp}\inv x x_{\alp}) \cdot f_{0} = 
	x_{\alp} \cdot \lpar (x_{\alp}\inv x x_{\alp}) \cdot f_{0} \rpar$. Since 
	$H$ is normal in $G$, we have $x_{\alp}\inv x x_{\alp} \in H$ and so 
	$(x_{\alp}\inv x x_{\alp}) \cdot f_{0} \in \orb_{0}$. Since $x \in H$ 
	is arbitrary, we conclude that
	\eq{inc}
	\orb_{\alp} \sset x_{\alp} \cdot \orb_{0}.
	\end{equation}
	It follows that $\ord{\orb_{\alp}} \leq \ord{\orb_{0}} = q^{-2} 
	\ord{\orb}$ (by \refp{char}) and so $\ord{\orb_{\alp}} < \ord{\orb}$. 
	By \refl{order}, we conclude that $\ord{\orb_{\alp}} = q^{-2} \ord{\orb}$ 
	and so $\ord{\orb_{\alp}} = \ord{\orb_{0}} = \ord{x_{\alp} \cdot 
	\orb_{0}}$. By \refe{inc}, we obtain $\orb_{\alp} = x_{\alp} \cdot 
	\orb_{0}$ as required.
	
	Now, suppose that $\bet \in \fq$ is such that $x_{\bet} \cdot f_{0} \in 
	\orb_{\alp}$. Then, there exists $x \in H$ such that $x_{\bet} \cdot 
	f_{0} = x \cdot (x_{\alp} \cdot f_{0}) = (x x_{\alp}) \cdot f_{0}$. It 
	follows that $x_{\bet}\inv x x_{\alp} \in C_{H}(f_{0})$. In particular, 
	$x_{\bet}\inv x x_{\bet} x_{\bet}\inv x_{\alp} \in H$ and so 
	$x_{\bet}\inv x_{\alp} \in H$ (because $H$ is normal in $G$ and $x 
	\in H$). Since $x_{\bet}\inv = 1 - \bet e + a$ for some $a \in 
	J^{2}$, we have $$x_{\bet}\inv x_{\alp} = 1 + (\alp - \bet) e - 
	(\alp \bet) e^{2} + \alp ea.$$ Since $J^{2} \sset U$, we conclude that 
	$(\alp - \bet) e \in U$ and this implies that $\alp = \bet$. It follows 
	that the $H$-orbits $\orb_{\alp} \sset \pi(\orb)$, for $\alp \in \fq$, 
	are all distinct. Hence, the union $\bcup_{\alp \in \fq} \orb_{\alp}$ 
	is disjoint and so $$\ord{\bcup_{\alp \in \fq} \orb_{\alp}} = \sum_{\alp 
	\in \fq} \ord{\orb_{\alp}} = q \cdot \ord{\orb_{0}} = q\inv \cdot 
	\ord{\orb} = \ord{\pi(\orb)}$$ (because $\orb$ is of type I, hence 
	$\ord{\orb} = q \ord{\pi(\orb)} = q^{2} \ord{\orb_{0}}$). Since 
	$\orb_{\alp} \sset \orb$ for all $\alp \in \fq$, we conclude that 
	$$\pi(\orb) = \bcup_{\alp \in \fq} \orb_{\alp}$$ as required.
	
	Finally, let $\Gam \sset G$ be a set of representatives of the cosets 
	of $H$ in $G$. Then $G$ is the disjont union $$G = \bcup_{x \in \Gam} 
	xH.$$ Since $\ord{G:H} = q$, we have $\ord{\Gam} = q$. Moreover, 
	for each $x \in \Gam$, there exists a unique $\alp \in \fq$ such 
	that $x \in x_{\alp} H$. It follows that $x \cdot \orb_{0} \sset 
	x_{\alp} \cdot \orb_{0}$ and, by order considerations, the equality must 
	occur.
		
	The proof is complete.
\end{proof}

Next, given an arbitrary $G$-orbit $\orb \in \Orb(G)$, we consider the 
restriction $(\phi_{\orb})_{H}$ of the class function $\phi_{\orb}$ to
to the maximal algebra subgroup $H$. For simplicity, we shall write 
$\phi = \phi_{\orb}$. We recall that, by definition (see \refe{class}), 
we have $$\phi(1+a) = \frac{1}{\sqrt{\ord{\orb}}} \sum_{f \in \orb} 
\psi_{f}(a)$$ for all $a \in J$.

Suppose that $\orb$ is of type I. Then, by \refp{dec}, we have 
$\pi(\orb) = \orb_{0}$. Let $a \in U$ be arbitrary and consider 
the class function $\phi_{\orb_{0}} \in \cf(H)$; for simplicity, we 
shall write $\phi_{0} = \phi_{\orb_{0}}$. We have $$\phi_{0}(1+a) = 
\frac{1}{\sqrt{\ord{\orb_{0}}}} \sum_{g \in \orb_{0}} \psi_{g}(a) = 
\frac{1}{\sqrt{\ord{\orb_{0}}}} \sum_{g \in \pi(\orb)} \psi_{g}(a).$$ 
On the other hand, since $\calL(f) \cap \orb = \lset f \rset$ for all 
$f \in \orb$ (by \refp{types}), the map $\pi$ determines naturally a 
bijection between the $G$-orbit $\orb \sset J^{*}$ and the $H$-orbit 
$\orb_{0} = \pi(\orb) \sset U^{*}$. Therefore $$\sum_{g \in \pi(\orb)} 
\psi_{g}(a) = \sum_{f \in \orb} \psi_{f}(a).$$ Since $\ord{\orb} = 
\ord{\orb_{0}}$ (by \refp{char}), we conclude that $$\phi(a) = 
\phi_{0}(a)$$ for all $a \in U$. It follows that $$\phi_{H} = 
\phi_{0}.$$

Now, suppose that $\orb$ is of type II. Then, by \refp{dec}, 
$\pi(\orb)$ is the disjoint union $$\pi(\orb) = \bcup_{\alp \in \fq} 
\orb_{\alp}$$ of $H$-orbits $\orb_{\alp}$ for $\alp \in \fq$. Let $a 
\in U$ be arbitrary. In this case, we have $$\sum_{g \in \pi(\orb)} 
\psi_{g}(a) = \sum_{\alp \in \fq} \sum_{g \in \orb_{\alp}} 
\psi_{g}(a) = \sum_{\alp \in \fq} \sqrt{\ord{\orb_{\alp}}} \cdot 
\phi_{\alp}(1+a)$$ where, for any $\alp \in \fq$, $\phi_{\alp}$ 
denotes the class function $\phi_{\orb_{\alp}} \in \cf(H)$. Since 
$\ord{\orb_{\alp}} =  q^{-2} \ord{\orb}$ for all $\alp \in \fq$ (see 
the proof of \refp{dec}), we conclude that $$\sum_{\alp \in \fq} 
\phi_{\alp}(1+a) = \frac{q}{\sqrt{\ord{\orb}}} \sum_{g \in \pi(\orb)} 
\psi_{g}(a).$$ On the other hand, we have $\calL(f) \sset \orb$ for all 
$f \in \orb$ (by \refp{types}). Hence, there exist elements $\seq{f}{r} 
\in \orb$ such that $\orb$ is the disjoint union $$\orb = \calL(f_{1}) 
\cup \ldots \cup \calL(f_{r}).$$ It follows that $$\phi(1+a) = 
\frac{1}{\sqrt{\ord{\orb}}} \sum_{f \in \orb} \psi_{f}(a) = 
\frac{1}{\sqrt{\ord{\orb}}} \sum_{i=1}^{r} \sum_{f \in \calL(f_{i})} 
\psi_{f}(a) = \frac{q}{\sqrt{\ord{\orb}}} \sum_{i=1}^{r} 
\psi_{f_{i}}(a)$$ (because $f(a) = f_{i}(a)$ for all $f \in 
\calL(f_{i})$). Finally, we clearly have $r = \ord{\pi(\orb)}$ and 
$\pi(\orb) = \lset \pi(f_{1}), \ldots, \pi(f_{r}) \rset$. Therefore, 
$$\sum_{i=1}^{r} \psi_{f_{i}}(a) = \sum_{g \in \pi(\orb)} 
\psi_{g}(a).$$ It follows that $$\phi_{H} = \sum_{\alp \in \fq} 
\phi_{\alpha}.$$ Finally, we note that, by \refp{dec}, for each $\alp 
\in \fq$, the $H$-orbit $\orb_{\alp}$ considered above may be chosen to 
be $x_{\alp} \cdot \orb_{0}$. Now, let $\alp \in \fq$ and $a \in U$ be 
arbitrary. Then, $$\phi_{\alp}(1+a) = \frac{1}{\sqrt{\ord{\orb_{\alp}}}} 
\sum_{g \in \orb_{\alp}} \psi_{g}(a) = \frac{1}{\sqrt{\ord{\orb_{0}}}} 
\sum_{f \in \orb_{0}} \psi_{x_{\alp} \cdot f}(a).$$ Let $f \in \orb_{0}$ 
be arbitrary. Then, by definition, $$\psi_{x_{\alp} \cdot f}(a) = 
\psi((x_{\alp} \cdot f)(a)) = \psi(f(x_{\alp}\inv a x_{\alp}) = 
\psi_{f}(x_{\alp}\inv a x_{\alp}).$$ This concludes the proof of the 
following.

\prop{rest}
	Let $H$ be a maximal algebra subgroup of $G$, let $U \sset J$ be such 
	that $H = 1+U$ and let $\apl{\pi}{J^{*}}{U^{*}}$ be the natural 
	projection. Let $\orb \in \Orb(G)$ be arbitrary and let $\phi$ 
	denote the class function $\phi_{\orb} \in \cf(G)$. Then, the 
	following statements hold:
	\begin{enumerate}
		\al{a}  If $\orb$ is of type I (with respect to $H$), then the 
		restriction $\phi_{H}$ of $\phi$ to $H$ is the class function 
		$\phi_{0} = \phi_{\orb_{0}} \in \cf(H)$ corresponding to the 
		$H$-orbit $\orb_{0} = \pi(\orb) \sset U^{*}$.
		\al{b}  Suppose that $\orb$ is of type II (with respect to $H$). 
		Let $\lset x_{\alp} \col \alp \in \fq \rset$ be a set of 
		representatives of the cosets of $H$ in $G$, let $\orb_{0} \in 
		\Orb(H)$ be an (arbitrary) $H$-orbit satisfying $\orb_{0} \sset 
		\pi(\orb)$ and, for each $\alp \in \fq$, let $\orb_{\alp} \in 
		\Orb(H)$ be the $H$-orbit $\orb_{\alp} = x_{\alp} \cdot \orb_{0}$. 
		Then, the restriction $\phi_{H}$ of $\phi$ to $H$ is the sum 
		$$\phi_{H} = \sum_{\alp \in \fq} \phi_{\alp}$$ of the (linearly 
		independent) class functions $\phi_{\alp} \in \cf(H)$, $\alp 
		\in \fq$, which correspond to the $H$-orbits $\orb_{\alp}$, 
		$\alp \in \fq$. Moreover, for each $\alp \in \fq$, $\phi_{\alp}$ is 
		the class function defined by $\phi_{\alp}(x) = \phi_{0}(x_{\alp}\inv 
		x x_{\alp})$ for all $x \in H$.
	\end{enumerate}
\eprop

In the next result, we use Frobenius reciprocity to obtain the 
decomposition of the class function ${\phi_{\orb_{0}}}^{G} \in \cf(G)$ 
induced from the class function $\phi_{\orb_{0}} \in \cf(H)$ which 
each associated with an arbitrary $H$-orbit $\orb_{0} \in \Orb(H)$.

\prop{ind}
	Let $H$ be a maximal algebra subgroup of $G$, let $U \sset J$ be such 
	that $H = 1+U$ and let $\apl{\pi}{J^{*}}{U^{*}}$ be the natural 
	projection. Let $\orb_{0} \in \Orb(H)$ be arbitrary and let $\orb 
	\in \Orb(G)$ be an arbitrary $G$-orbit satisfying $\orb_{0} \sset 
	\pi(\orb)$. Moreover, let $\phi_{0}$ denote the class function 
	$\phi_{\orb_{0}} \in \cf(H)$. Then, the following statements hold:
	\begin{enumerate}
		\al{a}  Suppose that $\orb$ is of type I (with respect to $H$). Let 
		$e \in J$ be such that $J = U \+ \fq e$ and let $e^{*} \in 
		U^{\perp}$ be such that $e^{*}(e) = 1$. Let $f \in \orb$ be 
		arbitrary and, for each $\alp \in \fq$, let $\orb(\alp) \in \Orb(G)$ 
		denote the $G$-orbit which contains the element $f + \alp e^{*} \in 
		J^{*}$. Then, the $G$-orbits $\orb(\alp)$, for $\alp \in \fq$, 
		are all distinct and the induced class function ${\phi_{0}}^{G}$ is 
		the sum $${\phi_{0}}^{G} = \sum_{\alp \in \fq} \phi_{\orb(\alp)}$$ of 
		the (linearly independent) class functions $\phi_{\orb(\alp)}$, for 
		$\alp \in \fq$.
		\al{b}  If $\orb$ is of type II (with respect to $H$), then 
		${\phi_{0}}^{G}$ is the class function $\phi = \phi_{\orb}$ which 
		corresponds to the $G$-orbit $\orb \in \Orb(G)$.
	\end{enumerate}
\eprop

\begin{proof}
	By \refp{ortho}, we have $${\phi_{0}}^{G} = \sum_{\orb\pr \in 
	\Orb(G)} \mu_{\orb\pr} \phi_{\orb\pr}$$ where $\mu_{\orb\pr} = 
	\lan {\phi_{0}}^{G}, \phi_{\orb\pr} \ran_{G}$ for all $\orb\pr \in 
	\Orb(G)$. Let $\orb\pr \in \Orb(G)$ be arbitrary. By Frobenius reciprocity, 
	we have $$\lan {\phi_{0}}^{G}, \phi_{\orb\pr} \ran_{G} = \lan \phi_{0}, 
	(\phi_{\orb\pr})_{H} \ran_{H}.$$ Therefore, by \refp{rest}, 
	$\mu_{\orb\pr} \neq 0$ if and only if $\orb_{0} \sset \pi(\orb\pr)$; 
	moreover, if this is the case, we have $\mu_{\orb\pr} = 1$. On the 
	other hand, let $f \in \orb$ be such that $\pi(f) \in \orb_{0}$; we 
	note that, in the case where $\orb$ is of type I (with respect to $H$), 
	we have $\pi(\orb) = \orb_{0}$ (by \refp{dec}) and so $\pi(f) \in 
	\orb_{0}$ for all $f \in \orb$. Let $\orb\pr \in \Orb(G)$ be such that 
	$\orb_{0} \sset \pi(\orb\pr)$. Then $\pi(f) = \pi(f\pr)$ for some $f\pr 
	\in \orb\pr$, hence $f\pr \in \calL(f)$. Since $U^{\perp} = \fq 
	e^{*}$, we have $\calL(f) = f + \fq e^{*}$ and so there exists $\alp \in 
	\fq$ such that $f\pr = f + \alp e^{*}$. Therefore, 
	$\orb\pr$ is the $G$-orbit which contains the element $f + \alp e^{*}$; 
	we denote this $G$-orbit by $\orb(\alp)$. 
	It follows that $${\phi_{0}}^{G} = \sum_{\alp \in \Gam} 
	\phi_{\orb(\alp)}$$ where $\Gam \sset \fq$ is such that the $G$-orbits 
	$\orb(\alp)$, for $\alp \in \Gam$, are all distinct.
	
	Suppose that $\orb$ is of type II. Then $\calL(f) \sset \orb$ (by 
	\refp{types}) and so $\orb(\alp) = \orb$ for all $\alp \in \fq$. It 
	follows that, in this case, $\ord{\Gam} = 1$ and so 
	$${\phi_{0}}^{G} = \phi$$ as required (in part (b)).
	
	On the other hand, suppose that $\orb$ is of type I. Let $\alp \in 
	\fq$ be arbitrary. Then, by \refp{types}, the $G$-orbit $\orb(\alp)$ 
	is also of type I; otherwise, $\calL(f) = \calL(f + \alp e^{*}) \sset 
	\orb(\alp)$. Therefore, by \refp{char}, $\ord{\orb(\alp)} = 
	\ord{\orb_{0}}$ and so $$q \sqrt{\ord{\orb_{0}}} = {\phi_{0}}^{G}(1) = 
	\sum_{\alp \in \Gam} \phi_{\orb(\alp)}(1) = \sum_{\alp \in \Gam} 
	\sqrt{\ord{\orb(\alp)}} = \ord{\Gam} \sqrt{\ord{\orb_{0}}}.$$ It 
	follows that $\ord{\Gam} = q$ and so $\Gam = \fq$. In particular, we 
	conclude that the $G$-orbits $\orb(\alpha)$, for $\alp \in \fq$, are 
	all distinct. Moreover, we obtain $${\phi_{0}}^{G} = \sum_{\alp \in 
	\fq} \phi_{\orb(\alp)}$$ as required (in part (a)).
	
	The proof is complete.
\end{proof}

As a consequence (of the proof) we deduce the following result.

\prop{inv}
	Let $H$ be a maximal algebra subgroup of $G$, let $U \sset J$ be such 
	that $H = 1+U$ and let $\apl{\pi}{J^{*}}{U^{*}}$ be the natural 
	projection. Let $\orb_{0} \in \Orb(H)$ be arbitrary and let $\orb 
	\in \Orb(G)$ be an arbitrary $G$-orbit satisfying $\orb_{0} \sset 
	\pi(\orb)$. Then, the following statements hold:
	\begin{enumerate}
		\al{a}  Suppose that $\orb$ is of type I (with respect to $H$). Let 
		$e \in J$ be such that $J = U \+ \fq e$ and let $e^{*} \in 
		U^{\perp}$ be such that $e^{*}(e) = 1$. Let $f \in \orb$ be 
		arbitrary and, for each $\alp \in \fq$, let $\orb(\alp) \in \Orb(G)$ 
		denote the $G$-orbit which contains the element $f + \alp e^{*} \in 
		J^{*}$. Then, the $G$-orbits $\orb(\alp)$, for $\alp \in \fq$, 
		are all distinct and the inverse image $\pi\inv(\orb_{0})$ decomposes 
		as the disjoint union $$\pi\inv(\orb_{0}) = \bcup_{\alp \in \fq} 
		\orb(\alp)$$ of the $G$-orbits $\orb(\alp) \in \Orb(G)$, for $\alp 
		\in \fq$.
		\al{b}  If $\orb$ is of type II (with respect to $H$), then 
		$\pi\inv(\orb_{0}) \sset \orb$. Moreover, this inclusion is proper. 
	\end{enumerate}
\eprop

\begin{proof}
	Suppose that $\orb$ is of type I. Let $\alp \in \fq$ be arbitrary. 
	Then, as we have seen in the proof of \refp{ind}, $\orb(\alp)$ is also 
	of type I and so $\pi(\orb(\alp)) = \orb_{0}$ (by \refp{dec} because 
	$\orb_{0} \sset \pi(\orb(\alp))$). It follows that $\orb(\alp) \sset 
	\pi\inv(\orb_{0})$. Conversely, suppose that $g \in J^{*}$ is such 
	that $\pi(g) \in \orb_{0}$. Since $\orb_{0} = \pi(\orb)$, there 
	exists $x \in G$ such that $\pi(g) = \pi(x \cdot f)$. It follows that 
	$g \in \calL(x \cdot f)$. Since $\calL(x \cdot f) = x \cdot 
	\calL(f)$ and since $\calL(f) = f + U^{\perp}$, we conclude that 
	$g = x \cdot \lpar f + \alp e^{*} \rpar$ (hence, $g \in \orb(\alp)$) 
	for some $\alp \in \fq$. This completes the proof of part (a).
	
	Now, suppose that $\orb$ is of type II. Let $g \in J^{*}$ be such 
	that $\pi(g) \in \orb_{0}$. Since $\orb_{0} \sset \pi(\orb)$, there 
	exists $f \in \orb$ such that $\pi(g) = \pi(f)$. Therefore, $g \in 
	\calL(f)$. Since $\orb$ is of type II, $\calL(f) \sset \orb$ (by 
	\refp{types}) and so $g \in \orb$. Since $g \in \pi\inv(\orb_{0})$ is 
	arbitrary, we conclude that $\pi\inv(\orb_{0}) \sset \orb$. To see 
	that this inclusion is proper, it is enough to choose $g \in \orb$ 
	such that $\pi(g) \in x_{\alp} \cdot \orb_{0}$ for some $\alp \in 
	\fq$, $\alp \neq 0$, where $\lset x_{\alp} \col \alp \in \fq \rset$ 
	is as in \refp{dec}.
\end{proof}

\sec{irred}{Irreducible characters}

The purpose of this section is the proof of the following result.

\th{irred}
	Let $G$ be an arbitrary $\fq$-algebra group. Then, for each $\orb 
	\in \Orb(G)$, the class function $\phi_{\orb}$ is an irreducible 
	character of $G$. Moreover, we have $\irr(G) = \lset \phi_{\orb} 
	\col \orb \in \Orb(G) \rset$.
\eth

By \refp{ortho}, it is enough to show that the class functions 
$\phi_{\orb} \in \cf(G)$, $\orb \in \Orb(G)$, are, in fact, 
characters. And, to see this, we proceed by induction on $\ord{G}$ (the 
result being clear if $\ord{G} = 1$). The key step is solved by the 
following lemma. As before, $G = 1+J$ where $J = J(A)$ is the 
Jacobson radical of a finite dimensional $\fq$-algebra $A$. 

\lem{irred}
	Let $H$ be a maximal algebra subgroup of $G$, let $U \sset J$ be 
	such that $H = 1 + U$ and let $\pi \col J^{*} \rar U^{*}$ be 
	the natural projection. Let $\orb \in \Orb(G)$ and let $\orb_{0} \in 
	\Orb(H)$ be any $H$-orbit with $\orb_{0} \sset \pi(\orb)$. Assume 
	that the class function $\phi_{\orb_{0}} \in \cf(H)$ is a character 
	of $H$. Then, the class function $\phi_{\orb} \in \cf(G)$ is a 
	character of $G$.
\elem

The proof of \refl{irred} relies on the following result (and on its 
corollary).

\lem{char}
	Let $H$ be a maximal algebra subgroup of $G$, let $\chi \in \irr(G)$ 
	and let $\tet \in \irr(H)$ be an irreducible constituent of $\chi_{H}$. 
	Then one (and only one) of the following two possibilities occurs:
	\begin{enumerate}
		\al{a}  $\chi_{H} = \tet$ is an irreducible character of $H$ and 
		$\tet^{G}$ has $q$ distinct irreducible constituents each one 
		occurring with multiplicity one. The irreducible constituents of 
		$\tet^{G}$ are the characters $\lam \chi$ for $\lam \in \irr(G/H)$ 
		(as usual, $\irr(G/H)$ is naturally identified with a subset of 
		$\irr(G)$).
		\al{b}  $\tet^{G} = \chi$ is an irreducible character of $G$ and 
		$\chi_{H}$ has $q$ distinct irreducible constituents each one 
		occurring with multiplicity one. 
	\end{enumerate}
\elem

\begin{proof}
	By \cite[Corollary~11.29]{isaacs2}, we know that $\chi(1)/\tet(1)$ 
	divides $\ord{G:H}$ (we recall that $H$ is a normal subgroup of 
	$G$). Since $G$ and $H$ are $\fq$-algebra subgroups, 
	\cite[Theorem~A]{isaacs1} asserts that $\chi(1)$ and $\tet(1)$ are 
	powers of $q$. Moreover, being a maximal algebra subgroup of $G$, we 
	know that $H$ has index $q$ in $G$. It follows that, either 
	$\chi(1) = \tet(1)$, or $\chi(1) = q \tet(1)$.
	
	Suppose that $\chi(1) = \tet(1)$. Then, we must have $\chi_{H} = 
	\tet$ and this is the situation of (a). The assertion concerning the 
	induced character is an easy application of a result of Gallagher (see 
	\cite[Corollary~6.17]{isaacs2}).
	
	On the other hand, suppose that $\chi(1) = q \tet(1)$. Then $\chi = 
	\tet^{G}$ (because $\tet^{G}(1) = q \tet(1)$ and because, by 
	Frobenius reciprocity, $\chi$ is an irreducible constituent of 
	$\tet^{G}$). By Clifford's Theorem (see, for example, 
	\cite[Theorem~6.2]{isaacs2}), we have $$\chi_{H} = e \sum_{i=1}^{t} 
	\tet_{i}$$ where $\tet = \tet_{1}, \tet_{2}, \ldots, \tet_{t}$ are 
	all the distinct conjugates of $\tet$ in $G$ and where $e = \lan 
	\chi_{H}, \tet \ran_{H}$. (The $G$-action on $\irr(H)$ is defined as 
	usual by $\tet^{x}(y) = \tet(x\inv y x)$ for all $\tet \in \irr(H)$, 
	all $x \in G$ and all $y \in H$.) By Frobenius reciprocity, we deduce 
	that $\lan \chi_{H}, \tet \ran_{H} = \lan \chi, \tet^{G} \ran_{G} = 
	1$ (because $\chi = \tet^{G}$). It follows that $\chi(1) = 
	\sum_{i=1}^{t} \tet_{i}(1) = t \tet(1)$ (because $\tet_{i}(1) = 
	\tet(1)$ for all $1 \leq i \leq t$) and so $t = q$.
	
	The proof is complete.
\end{proof}

The following easy consequence will also be very useful.

\cor{frob}
	Let $H$ be a maximal algebra subgroup of $G$ and let $\chi, \chi\pr \in 
	\irr(G)$. Then, the following hold:
	\begin{enumerate}
		\al{a}  If $\chi_{H}$ is irreducible, then $\lan \chi_{H}, \chi\pr_{H} 
		\ran_{H} \neq 0$ if and only if $\chi\pr = \lam \chi$ for some $\lam 
		\in \irr(G/H)$.
		\al{b}  If $\chi_{H}$ is reducible, then $\lan \chi_{H}, \chi\pr_{H} 
		\ran_{H} \neq 0$ if and only if $\chi\pr = \chi$. Moreover, we have 
		$\lan \chi_{H}, \chi_{H} \ran_{H} = q$.
	\end{enumerate}
\ecor

\begin{proof}
	Suppose that $\chi_{H} = \tet \in \irr(H)$. By Frobenius reciprocity, 
	we have $\lan \tet, \chi\pr_{H} \ran_{H} = \lan \tet^{G}, \chi\pr 
	\ran_{G}$ and so $\lan \chi_{H}, \chi\pr_{H} \ran_{H} \neq 0$ if and only 
	if $\chi\pr = \lam \chi$ for some $\lam \in \irr(G/H)$ (by part (a) of 
	\refl{char}).
	
	On the other hand, suppose that $\chi_{H}$ is reducible and let 
	$\chi\pr \in \irr(G)$ be such that $\lan \chi_{H}, \chi\pr_{H} \ran_{H} 
	\neq 0$. Then, there exists $\tet \in \irr(H)$ such that $\lan 
	\tet, \chi_{H} \ran_{H} \neq 0$ and $\lan \tet, \chi\pr_{H} \ran_{H} 
	\neq 0$. By part (b) of \refl{char}, we have $\chi = \tet^{G}$ and 
	so, using Frobenius reciprocity, we deduce that $\lan \chi, \chi\pr 
	\ran_{G} = \lan \tet^{G}, \chi\pr \ran_{G} = \lan \tet, \chi\pr_{H} 
	\ran_{H} \neq 0$. Since $\chi, \chi\pr \in \irr(G)$, we conclude 
	that $\chi = \chi\pr$. Finally, by part (b) of \refl{char}, it is clear 
	that $\lan \chi_{H}, \chi_{H} \ran_{H} = q$.
	
	The proof is complete.
\end{proof}

We are now able to prove \refl{irred}.

\renewcommand{\proofname}{Proof of \refl{irred}}

\begin{proof}
	For simplicity, we write $\phi = \phi_{\orb}$ and $\phi_{0} = 
	\phi_{\orb_{0}}$. If $\orb$ is of type II with respect to $H$, then 
	$\phi = (\phi_{0})^{G}$ (by \refp{ind}) and the result is clear. On 
	the other hand, suppose that $\orb$ is of type I with respect to $H$. 
	Then, by \refp{rest}, $\phi_{0} = \phi_{H}$ and $\lan {\phi_{0}}^{G}, 
	{\phi_{0}}^{G} \ran_{G} = q$. By \refp{ortho}, we have $$\phi = 
	\sum_{\chi \in \irr(G)} \mu_{\chi} \chi$$ where $\mu_{\chi} \in \CC$ for 
	all $\chi \in \irr(G)$. Let $\calI = \lset \chi \in \irr(G) \col 
	\mu_{\chi} \neq 0 \rset$ be the support of $\phi$. Since $\lan \phi, \phi 
	\ran_{G} = 1$, we have
	\eq{mu}
		\sum_{\chi \in \calI} \ord{\mu_{\chi}}^{2} = 1.
	\end{equation}
	Now, we consider the restriction
	\eq{rest}
		\phi_{0} = \phi_{H} = \sum_{\chi \in \calI} \mu_{\chi} \chi_{H}
	\end{equation}
	of $\phi$ to $H$.
	
	We claim that $\chi_{H} \in \irr(H)$ for all $\chi \in \calI$. To see 
	this, suppose that $\chi_{H}$ is reducible for some $\chi \in \calI$. 
	Let $\tet \in \irr(H)$ be an irreducible constituent of $\chi_{H}$. Then, 
	$\tet$ occurs in the sum of the right hand side of \refe{rest} with 
	coefficient $\mu_{\chi}$; in fact, $\lan \tet, \chi_{H} \ran_{H} = 
	1$ (by \refl{char}) and $\lan \tet, \chi\pr_{H} \ran_{H} = 0$ for 
	all $\chi\pr \in \irr(G)$ with $\chi\pr \neq \chi$ (otherwise, $\lan 
	\chi_{H}, \chi\pr_{H} \ran_{H} \neq 0$ and this is in contradiction 
	with \refc{frob}). It follows that \refe{rest} has the form
	\eq{rest2}
		\phi_{0} = \mu_{\chi} \tet_{1} + \cdots + \mu_{\chi} \tet_{q} + \xi
	\end{equation}
	where $\tet_{1} = \tet, \tet_{2}, \ldots, \tet_{q}$ are the $q$ distinct 
	irreducible constituents of $\chi_{H}$ (see \refl{char}) and where $\xi 
	\in \cf(H)$ is a $\CC$-linear combination of $\irr(H) \setminus \lset 
	\seq{\tet}{q} \rset$. Now, since $\phi_{0}$ is a character of $H$ (by 
	assumption) and since $\lan \phi_{0}, \phi_{0} \ran_{H} = 1$ (by 
	\refp{ortho}), we have $\phi_{0} \in \irr(H)$. Since $\irr(H)$ is a 
	$\CC$-basis of $\cf(H)$, the equality \refe{rest2} implies that 
	$\mu_{\chi} = 0$ and this is in contradiction with $\chi \in \calI$. This 
	completes the proof of our claim, \ie $\chi_{H} \in \irr(H)$ for all 
	$\chi \in \calI$.
	
	Now, for each $\tet \in \irr(H)$, let $\calI_{\tet} = \lset \chi \in 
	\calI \col \chi_{H} = \tet \rset$ and let $$\mu_{\tet} = \sum_{\chi \in 
	\calI_{\tet}} \mu_{\chi}.$$ Then $\calI$ is the disjoint union $$\calI = 
	\bcup_{\tet \in \irr(H)} \calI_{\tet}$$ and so $$\phi_{0} = \sum_{\tet 
	\in \irr(H)} \mu_{\tet} \tet.$$ Since $\phi_{0} \in \irr(H)$, we conclude 
	that $\mu_{\tet} = \del_{\tet, \phi_{0}}$ for all $\tet \in \irr(H)$. 
	Hence, $$\sum_{\chi \in \calI} \mu_{\chi} = \sum_{\tet \in \irr(H)} 
	\mu_{\tet} = 1.$$ Using \refe{mu}, we easily deduce that there exists 
	a unique $\chi \in \irr(G)$ with $\mu_{\chi} \neq 0$ and, in fact, 
	$\mu_{\chi} = 1$.
	
	The proof is complete.
\end{proof} 

\renewcommand{\proofname}{Proof}

\sec{proof}{Proof of \reft{main}}

In this section, we prove \reft{main}. As before, let $A$ be a 
finite dimensional $\fq$-algebra, let $J = J(A)$ be the Jacobson radical 
of $A$ and let $G = 1+J$ be the $\fq$-algebra group defined by $J$. 

We consider the chain $J \Sset J^{2} \Sset J^{3} \Sset \ldots$ of ideals 
of $A$. Since $J$ is nilpotent, there exists the smallest integer $m$ with 
$J^{m} \neq \lset 0 \rset$. Moreover, we may refine the chain $$J \Sset 
J^{2} \Sset \ldots \Sset J^{m} \Sset \lset 0 \rset$$ to obtain a (maximal) 
chain $$\lset 0 \rset = U_{0} \sset U_{1} \sset \ldots \sset U_{n} = 
J$$ of ideals of $A$ satisfying $$\dim U_{i+1} = \dim U_{i}+1$$ for 
all $0 \leq i < n$. Let $f \in J^{*}$ be arbitrary and, for each $0 \leq 
i \leq n$, let $$R_{i} = \lset a \in U_{i} \col f([ab]) = 0 \all b \in 
U_{i} \rset.$$ Finally, let $$U = R_{1} + \cdots + R_{n}.$$

It is clear that $U$ is an $\fq$-subspace of $J$. Now, let $a,b \in U$ 
and suppose that $a \in R_{i}$ and $b \in R_{j}$ for some $1 \leq i, j 
\leq n$ with $i \leq j$. We claim that $ab \in R_{i}$. To see this, let 
$c \in U_{i}$ be arbitrary. Then, $$f([ab,c]) = f([a,bc]) + f([b,ca]).$$ 
Since $U_{i}$ is an ideal of $A$, we have $bc \in U_{i}$, hence 
$f([a,bc]) = 0$ (because $a \in R_{i}$). On the other hand, we have $ca 
\in U_{i}$. Since $U_{i} \sset U_{j}$ (because $i \leq j$) and and 
since $b \in R_{j}$, we conclude that $f([b,ca]) = 0$. Thus $$f([ab,c]) = 
f([a,bc]) + f([b,ca]) = 0$$ and so $ab \in R_{i}$ (because $c \in U_{i}$ 
is arbitrary). It follows that $U$ is a multiplicatively closed 
$\fq$-subspace of $J$. Hence, $H = 1+U$ is an algebra subgroup of $G$.
Moreover, a similar argument shows that $$f([ab]) = 0$$ for all $a,b \in 
U$. This means that $U$ is an {\it $f$-isotropic} $\fq$-subspace of 
$J$ (\ie $U$ is isotropic with respect to the skew-symmetric bilinear form 
$B_{f}$ which was defined in \refs{class}.) Next, we claim that $U$ is 
a maximal $f$-isotropic $\fq$-subspace of $J$. By Witt's Theorem (see, 
for example, \cite[Theorems~3.10~and~3.11]{artin}), it is enough to 
prove that
\eq{mdim}
	\dim U = \frac{1}{2} \lpar \dim J + \dim \rad(f) \rpar.
\end{equation}
To see this, we proceed by induction on $\dim J$. If $\dim J = 1$, then 
$U = J = \rad(f)$ and the claim is trivial. Now, suppose that $\dim J > 
1$ and consider the ideal $U_{n-1}$ of $J$. Let $$U\pr = R_{1} + \cdots + 
R_{n-1}.$$ By induction, we have $$\dim U\pr = \frac{1}{2} \lpar \dim 
U_{n-1} + \dim \rad(f\pr) \rpar$$ where $f\pr$ is the restriction of $f$ 
to $U_{n-1}$. Using \refp{rad} and \refp{char}, we conclude that, either 
$\dim \rad(f) = \dim \rad(f\pr) - 1$, or $\dim \rad(f) = \dim 
\rad(f\pr) + 1$. In the first case, we deduce that $$\dim U\pr = 
\frac{1}{2} \lpar \dim J - 1 + \dim \rad(f) + 1 \rpar = \frac{1}{2} \lpar 
\dim J + \dim \rad(f) \rpar.$$ Therefore, $U\pr$ is a maximal $f$-isotropic 
$\fq$-subspace of $J$. Since $U\pr \sset U$, we conclude that $U\pr = U$ 
and \refe{mdim} follows in this case. On the other hand, suppose that 
$\dim \rad(f) = \dim \rad(f\pr) + 1$. Then, $$\dim U\pr = \frac{1}{2} \lpar 
\dim J - 1 + \dim \rad(f) - 1 \rpar = \frac{1}{2} \lpar \dim J + \dim 
\rad(f) \rpar - 1.$$ Since $U\pr \sset U$, we have $\dim U\pr \leq \dim 
U$. If $\dim U\pr = \dim U$, then $U\pr = U$ and so $\rad(f) \sset U\pr 
\sset U_{n-1}$. If this were the case, we should have $\rad(f) \sset 
\rad(f\pr)$ and so $\dim \rad(f) \leq \dim \rad(f\pr)$, a 
contradiction. It follows that $$\dim U\pr < \dim U.$$ Since $U$ is an 
$f$-isotropic $\fq$-subspace of $J$, we deduce that $$\dim U\pr \leq \dim 
U - 1 \leq \frac{1}{2} \lpar \dim J + \dim \rad(f) \rpar - 1 = \dim U\pr.$$ 
The proof of \refe{mdim} is complete.

Given an arbitrary element $f \in J^{*}$, we will say that a multiplicatively 
closed $\fq$-subspace $U$ of $J$ is an {\bf $f$-polarization} if $U$ is 
a maximal $f$-isotropic $\fq$-subspace of $J$. Hence, we have finished 
the proof of the following.

\prop{polar}
	Let $f \in J^{*}$ be arbitrary. Then, there exists a $f$-polariza\-tion 
	$U \sset J$.
\eprop

Now, it is easy to conclude the proof of \reft{main}. 

\renewcommand{\proofname}{Proof of \reft{main}}

\begin{proof}
	Let $\chi$ be an (arbitrary) irreducible character of $G$. Then, by 
	\reft{irred}, $\chi = \phi_{\orb}$ for some $G$-orbit $\orb \in 
	\Orb(G)$. Let $f \in \orb$ be arbitrary and let $U \sset J$ be an 
	$f$-polarization. Then, $H = 1+U$ is an algebra subgroup of $G$. Let 
	$f_{0} \in U^{*}$ be the restriction of $f$ to $U$. Since $U$ is 
	$f$-isotropic, we have $\rad(f_{0}) = U$, hence $C_{H}(f_{0}) = H$ (by 
	\refp{rad}). It follows that $\orb_{0} = \lset f_{0} \rset$ is a 
	single $H$-orbit on $U^{*}$ (\ie an element of $\Orb(H)$). We denote 
	by $\lam_{f}$ the class function $\phi_{\orb_{0}}$ of $H$; by 
	definition, $\lam_{f} \col H \rar \CC$ is defined by 
	$$\lam_{f}(1+a) = \psi_{f}(a) = \psi(f(a))$$ for all $a \in U$. By 
	\reft{irred}, we know that $\lam_{f}$ is an irreducible 
	character of $H$. Moreover, $\lam_{f}(1) = \sqrt{\ord{\orb_{0}}} = 1$, 
	\ie $\lam_{f}$ is a linear character of $H$. To conclude the proof, 
	we claim that
	\eq{induced}
		\phi_{\orb} = {\lam_{f}}^{G}.
	\end{equation}
	To see this, we evaluate the Frobenius product $\lan \phi_{\orb}, 
	{\lam_{f}}^{G} \ran_{G}$. Using Frobenius reciprocity (and the 
	definition of $\phi_{\orb}$), we deduce that
	\begin{eqnarray*}
		\lan \phi_{\orb}, {\lam_{f}}^{G} \ran_{G} & = & \lan 
		(\phi_{\orb})_{H}, \lam_{f} \ran_{H}  \\
		 & = & \frac{1}{\ord{H}} \sum_{x \in H} \phi_{\orb}(x) 
		 \ovl{\lam_{f}(x)}  \\
		 & = & \frac{1}{\ord{U}} \sum_{a \in U} \lpar 
		 \frac{1}{\sqrt{\ord{\orb}}} \sum_{g \in \orb} \psi_{g}(a) \rpar 
		 \ovl{\psi_{f}(a)}  \\
		 & = & \frac{1}{\sqrt{\ord{\orb}}} \sum_{g \in \orb} \lpar 
		 \frac{1}{\ord{U}} \sum_{a \in U} \psi_{g}(a) \ovl{\psi_{f}(a)} \rpar  \\
		 & = & \frac{1}{\sqrt{\ord{\orb}}} \lan \psi_{g}, \psi_{f} \ran_{U^{+}}.
	\end{eqnarray*}
	By \refe{irrJ}, given any $g \in J^{*}$, we have $\lan \psi_{g}, 
	\psi_{f} \ran_{U^{+}} \neq 0$ if and only if $f(a) = g(a)$ for all $a \in 
	U$; in other words, $\lan \psi_{g}, \psi_{f} \ran_{U^{+}} \neq 0$ if and 
	only if $g \in f+U^{\perp}$. Since $\psi_{g}$ is linear for all $g \in 
	J^{*}$, we conclude that
	\eq{frobenius}
		\lan \phi_{\orb}, {\lam_{f}}^{G} \ran_{G} = \frac{\ord{(f+U^{\perp}) 
		\cap \orb}}{\sqrt{\ord{\orb}}}.
	\end{equation}
	It follows that $\lan \phi_{\orb}, {\lam_{f}}^{G} \ran_{G} \neq 0$ 
	(because $f \in (f+U^{\perp}) \cap \orb$), hence $\phi_{\orb}$ is an 
	irreducible constituent of ${\lam_{f}}^{G}$. On the other hand, we 
	have $\phi_{\orb}(1) = \sqrt{\ord{\orb}}$ and
	\begin{eqnarray*}
		{\lam_{f}}^{G}(1) & = & \ord{G:H} \lam_{f}(1) = q^{\dim J - \dim U}  \\
		 & = & \sqrt{q^{\dim J - \dim \rad(f)}} = \sqrt{\ord{G:C_{G}(f)}} = 
		 \sqrt{\ord{\orb}}
	\end{eqnarray*}
	(using \refe{mdim} and \refp{rad}). The claim \refe{induced} follows and 
	the proof of \reft{main} is complete.
\end{proof}

\renewcommand{\proofname}{Proof}

\end{document}